\documentclass{gen-j-l}

\usepackage{amssymb, amsfonts}
\usepackage{amsmath}
\usepackage{xcolor}
\usepackage{mathtools}
\usepackage{enumerate}
\usepackage{stmaryrd}
\usepackage{amsthm}
\usepackage[left=3.25cm,right=3.25cm,top=2.5cm,bottom=2.5cm]{geometry}
\usepackage{relsize}
\usepackage{comment}

\newtheorem{theorem}{Theorem}[section]
\newtheorem{lemma}[theorem]{Lemma}
\newtheorem{cor}[theorem]{Corollary}
\newtheorem*{lemma*}{Lemma}

\theoremstyle{definition}

\newtheorem{example}[theorem]{Example}

\theoremstyle{remark}
\newtheorem{remark}[theorem]{Remark}

\numberwithin{equation}{section}

\newcommand{\R}{\mathbb{R}}  
\DeclareMathOperator{\totvar}{Tot.Var.} 

\begin{document}

\title[Error estimates of FV and RKDG methods for hyperbolic systems in 1D]{Conditional a priori error estimates of finite volume and Runge-Kutta discontinuous Galerkin methods with abstract limiting for hyperbolic systems of conservation laws in 1D}


\author{Fabio Leotta}
\address{}
\curraddr{}
\email{leotta@mathematik.tu-darmstadt.de}
\thanks{}

\subjclass[2020]{Primary }

\date{2024}

\dedicatory{}


\begin{abstract}
	We derive conditional a priori error estimates of a wide class of finite volume and Runge-Kutta discontinuous Galerkin methods with abstract limiting for hyperbolic systems of conservation laws in 1D via the verification of weak consistency and entropy stability, as recently proposed by Bressan et al.~\cite{BressanChiriShen21}. Convergence in $L^\infty L^1$ with rate $h^{1/3}$ is obtained under a time step restriction $\tau\leq ch$, provided the following conditions hold: the exact solution is piecewise Lipschitz continuous, its (finitely many and isolated) shock curves can be traced with precision $h^{2/3}$ and, outside of these shock tracing tubular neighborhoods the numerical solution -- assumed to be uniformly small in BV -- has oscillation strength $h$ across each mesh cell and cell boundary.
\end{abstract}

\maketitle

\tableofcontents

\section{Introduction}\label{sec:1}
We consider a strictly hyperbolic system of $n$ conservation laws in one space dimension,
\begin{align}\label{eq:cl}
	\begin{cases}
		\partial_t u + \partial_x f(u) = 0,& (t,x)\in(0,T]\times\R\\
		 u(0,x) = u_0(x),& x\in\R
	\end{cases},
\end{align}
with each characteristic field being either linearly degenerate or genuinely nonlinear. Here, $u:[0,T]\times\R\to\mathcal{O}\subset\R^n$ denotes the conserved quantity, $f\in C^2_b(\mathcal{O};\R^n)$ and $u_0\in L^\infty(\R;\mathcal{O})$ is a compactly supported initial datum that is piecewise Lipschitz continuous. The system (\ref{eq:cl}) shall be endowed with a strictly convex entropy $\eta$ and entropy flux $q$ such that
\begin{equation}\label{eq:entropy compatibility}
	\nabla\eta Df=\nabla q.
\end{equation} 
\begin{remark}
For future reference, we note here that the compatibility condition, (\ref{eq:entropy compatibility}), yields, after differentiation, the symmetry
\begin{equation}\label{eq:symmetry}
	\nabla^2\eta Df = (Df)^\intercal\nabla^2\eta.
\end{equation}	
\end{remark}

We will operate from within the well-posedness theory of small BV solutions, that is, the system (\ref{eq:cl}) generates a semigroup of entropy weak solutions $\mathcal{S}:\mathcal{D}\times[0,T]\to\mathcal{D}$, where
\begin{equation}\label{eq:D}
	\mathcal{D}:=\{\overline{u}\in L^1(\R;\R^n)~|~\totvar(\overline{u})\leq\delta_0\},
\end{equation}
for a certain $\delta_0>0$ and we write $\mathcal{S}(u,t)\hat{=}\mathcal{S}_t(u)$. In particular, the semigroup is Lipschitz continuous in the sense that there exist constants $C_0,L_0>0$ such that
\begin{align}
	\|\mathcal{S}_t(u)-\mathcal{S}_s(u)\|_{L^1(\R)}\leq& C_0|t-s|\totvar(u),\\
	\|\mathcal{S}_t(u)-\mathcal{S}_t(v)\|_{L^1(\R)}\leq& L_0\|u-v\|_{L^1(\R)},
\end{align}
for all $u,v\in\mathcal{D}$ and $s,t\in[0,T]$; cf.~\cite{BressanBook} or \cite{BianchiniBressan05} for an existence proof via front tracking or vanishing viscosity approximation, respectively. Lastly, we assume that the exact solution to the Cauchy problem (\ref{eq:cl}) vanishes outside a compact set of space-time $[0,T]\times\Omega$ where $\Omega\subset\R$ denotes a bounded interval.

Contrary to the scalar setting $n=1$, the question of convergence of fully discrete approximations to (\ref{eq:cl}) remains largely open. Generally speaking, the lack of convergence results can be attributed to the fact that it is notoriously difficult to prove strong BV-stability properties of finite volume schemes (or rather, any numerical method that is not a wave front tracking method) for (\ref{eq:cl}), although numerical evidence indeed suggests that this holds true in many cases; cf.~e.g.~the finite volume methods studied in the monograph of Toro \cite{Eleuterio09}. However, the question of BV-stability remains a fundamental obstacle to the proof of \textit{unconditional convergence} results since it is unclear under which conditions the resonance effect studied by Bressan and coworkers in \cite{BressanJenssenBaiti06} producing oscillations in a Godunov approximation can be prevented. Accordingly, known convergence results apply only to special systems, which we briefly review in the following.

For certain $2\times 2$ systems -- sometimes called straight line or Temple systems -- LeVeque and Temple \cite{LeVequeTemple85} proved BV stability and thus convergence of a subsequence for Godunov's method. This result has since been generalized to $n\times n$ straight line systems by Bressan and Jenssen \cite{BressanJenssen00}. For the isentropic Euler system, convergence of a subsequence can be proven for the Lax-Friedrichs or Godunov method via compensated compactness as done by Ding et al.~\cite{DingChenLuo89}. However, in all of the above cases, no convergence rate is retrieved. For linear systems, (\ref{eq:cl}) can be decoupled into $n$ independent scalar conservation laws and thus the scalar theory for error estimates easily carries over; cf.~\cite{CockburnLinShu89}. 

A somewhat different approach can be taken when generalized solution concepts are considered, especially interesting in the multi-dimensional setting where non-uniqueness of weak entropy solutions has been proven (e.g.~via convex integration, cf.~\cite{Markfelder21}) and is popularly attributed to turbulence. In the monograph of Feireisl, Luk\'acov\'a-Medvidov\'a and coworkers \cite{MariaEduardBook}, (weak) convergence of finite volume methods for the (2D and 3D) Euler system is proven via dissipative measure-valued solutions. If the exact solution of the Euler system is smooth enough, the convergence is indeed strong and an a priori error estimate can be given via the relative energy method.

Thus, if the exact solution to (\ref{eq:cl}) is not smooth, e.g.~in the presence of shocks, all of the aforementioned convergence theory is not amenable to deriving a priori error estimates, even for the special systems considered there. 

Recently, Giesselmann and Sikstel \cite{GiesselmannSikstel25} have provided a posteriori error estimates for first-order finite-volume approximations to (\ref{eq:cl}) via stability results obtained by Bressan et al. \cite{BressanChiriShen21}. In \cite{BressanChiriShen21}, acknowledging the fundamental obstructions to the derivation of unconditional a priori error estimates as described in the preceding paragraphs, an error estimation framework for abstract numerical methods for hyperbolic systems in 1D is proposed: Assuming that the numerical solution - at the discretization level in question - is $\text{(a) small in BV}$, and (b) does not oscillate too much outside a finite number of narrow strips in $[0,T]\times \R$, they prove that if the numerical method is consistent and entropy stable in a $(W^{1,\infty})^\ast$-sense, the numerical error at the current discretization level can be bounded by a computable quantity that scales with the mesh parameters times an unknown uniform constant originating from the stability theory. However, if (a) and (b) hold uniformly in the discretization parameters, the scaling of the error bound w.r.t.~the mesh parameters is definitive and one thus obtains an a priori error estimate. In this sense, one can speak of a \textit{conditional convergence} result for numerical schemes that yields an explicit convergence rate, akin to a quantitative version of the well-known Lax-Wendroff theorem. More specifically, convergence in $L^\infty L^1$ with rate $h^{1/3}$ is obtained under a time step restriction $\tau\leq ch$, provided the following conditions hold: the exact solution is piecewise Lipschitz continuous, its (finitely many and isolated) shock curves can be traced with precision $h^{2/3}$ and, outside of these shock tracing tubular neighborhoods the numerical solution -- apart from being uniformly small in BV and weakly consistent and entropy stable -- has oscillation strength $h$ across each mesh cell and cell boundary; cf.~Section \ref{sec:2} for a precise statement. 

Regarding the necessity to estimate locations of large jumps in the numerical solution, let us remark that although such information can be extracted from finite volume schemes only heuristically, jump detection strategies have long been used successfully in practice for mesh refinement; cf.~e.g.~\cite{DahmenGmSiggi01} and \cite{SiggiAleksey22} for discontinuity indicators via multi resolution analysis for finite volume and discontinuous Galerkin methods, respectively.

In this paper, we show that the abstract \textit{conditional a priori error estimates} derived from \cite{BressanChiriShen21} applies to a broad class of finite volume (FV) and Runge-Kutta discontinuous Galerkin (RKDG) methods with abstract limiting. This is done via the verification of weak consistency and entropy stability (cf.~(\ref{eq:LC})-(\ref{eq:WES}) for a definition of these notions) of the schemes.  In \cite{BressanChiriShen21}, weak consistency and entropy stability of the Lax-Friedrichs and Godunov method is proven by using the fact that these schemes are essentially defined by iteratively averaging and evolving these averages exactly according to (\ref{eq:cl}). The situation for general FV and RKDG schemes is fundamentally different, and thus so is our strategy which in particular hinges on monotonicity properties of the numerical flux; examples include the local Lax-Friedrichs and HLL flux with entropy fix.

The paper is structured as follows. First, we introduce the error estimation framework from \cite{BressanChiriShen21} in Section~\ref{sec:2} to formulate the convergence theorem that we want to apply. In Section \ref{sec:fvmethod}, we describe the class of FV methods under consideration, which we analyze in Section \ref{sec:fvanalysis}. Building on that, we consider forward Euler discontinuous Galerkin methods (without limiting) in Section \ref{sec:without limiting} as a preliminary step before turning to RKDG methods with abstract limiting in Section \ref{sec:with limiting}. Corollary~\ref{cor:fv} and Corollary~\ref{cor:DG} constitute the main results of this paper.

\section{The error estimation framework}\label{sec:2}
In this section, we briefly articulate the desired properties of an approximate solution that will then yield an error estimate as proven in \cite{BressanChiriShen21}. 

Let us introduce the equidistant space grid $x_{j+1/2}=jh,~j\in\mathbb{Z}$, with $h>0$ the mesh width, and let us label the cells as $I_j=(x_{j-1/2},x_{j+1/2})$ with midpoint $x_j=(x_{j-1/2}+x_{j+1/2})/2$. Furthermore, the time grid is given by $t^n=n\tau$, where $n\in\{0,1,\ldots,N\}$ and $N\tau=T$.

Note that here and in the following we will use the notation $x\lesssim y$ if there exists a constant $c>0$ independent of $\tau,h$ and possibly a test function $\varphi$ such that $x\leq c y$.

In this paper, we will prove that the considered numerical schemes yield approximate\footnote{In contrast to the relative energy method, no reconstruction procedure is necessary to render e.g.~finite volume solutions Lipschitz continuous; it suffices to work with the piecewise constant (in space and time) FV solution directly.} solutions $u~\widehat{=}~u_{\tau,h}$ to (\ref{eq:cl}) that fulfill the following three inequalities:
\begin{itemize}
	\item [(\textbf{LC})] The mapping $t\mapsto u(t,\cdot)\in L^1(\R;\R^n)$ is Lipschitz continuous when restricted to the time grid, i.e., for every pair of time levels $n<m$, one has
	\begin{equation}\label{eq:LC}
		\|u(t^m,\cdot)-u(t^n,\cdot)\|_{L^1(\R)}\lesssim (t^m-t^n)\sup\limits_{t\in[t^n,t^m]} \totvar\{u(t,\cdot)\}.
	\end{equation}
	\item [(\textbf{WC})] The approximate solution is weakly consistent, i.e., for every pair of time levels $n<m$, and every $\varphi\in C^1_c(\R^2)$, one has
	\begin{align}
		&\left|\int_{t^n}^{t^m}\int_\R u\varphi_t + f(u)\varphi_x~dx~dt - \int_\R u(t^m,x)\varphi(t^m,x) - u(t^n,x)\varphi(t^n,x)~dx\right|\nonumber\\
		\lesssim&~ h\|\varphi\|_{W^{1,\infty}(\R^2)}(t^m-t^n)\sup\limits_{t\in[t^n,t^m]} \totvar\{u(t,\cdot)\}\label{eq:WC}.
	\end{align}
	
	\item [(\textbf{WES})] The approximate solution is weakly entropy stable, i.e., for every pair of time levels $n<m$, and every $\varphi\in C^1_c(\R^2)$ with $\varphi\geq 0$, one has
	\begin{align}
		&\int_{t^n}^{t^m}\int_\R \eta(u)\varphi_t + q(u)\varphi_x~dx~dt\nonumber\\
		&-\int_\R \eta(u(t^m,x))\varphi(t^m,x) - \eta(u(t^n,x))\varphi(t^n,x)~dx\nonumber\\
		\gtrsim& -h\|\varphi\|_{W^{1,\infty}(\R^2)}(t^m-t^n)\sup\limits_{t\in[t^n,t^m]} \totvar\{u(t,\cdot)\}\label{eq:WES}.
	\end{align}
\end{itemize}

Furthermore, for the conditional a priori error estimates to hold, we assume the following $BV$-control hypothesis:
\begin{itemize}
	\item [(\textbf{BV})] The numerical solution is $BV$-stable in the sense that,
	\begin{equation}
		\sup\limits_{t\in[0,T]}\totvar(u(t,\cdot))\leq \delta_0,
	\end{equation} 
	holds uniformly in $(\tau,h)$, with $\delta_0$ the constant in \ref{eq:D}. Additionally, there exists finitely many Lipschitz curves $(t,\gamma(t))$ away from which the numerical solution has local oscillation strength $h$; more precisely, for all such curves $\gamma$ and $t\in[0,T]$, $j\in\mathbb{Z}$ such that \mbox{$|x_j-\gamma(t)|\gtrsim h^{2/3}$}, there holds 
	\begin{equation}
		\|u_x\|_{L^1(I_j)}+|\llbracket u\rrbracket|_{j-1/2}+|\llbracket u\rrbracket|_{j+1/2}\lesssim h,
	\end{equation}
	for all $\tau, h>0$ small enough.
\end{itemize}

\begin{remark}
	In the following sections, we will prove (\textbf{LC}), (\textbf{WC}) and (\textbf{WES}) for a wide class of numerical schemes. The BV-control property (\textbf{BV}), on the other hand, should be understood as our working hypothesis, as it is notoriously difficult to prove for finite volume type schemes approximating general hyperbolic systems in 1D, although well-designed numerical schemes for (\ref{eq:cl}) indeed appear to behave in a BV-stable fashion in many cases; cf.~the multitude of examples presented in \cite{Eleuterio09}.  However, as mentioned in the introduction, the instance of BV-instability for Godunov's method studied in \cite{BressanJenssenBaiti06} poses as a fundamental counterexample. It remains unclear under which conditions the resonance effect causing oscillation growth in \cite{BressanJenssenBaiti06} can be prevented.
\end{remark}

\begin{theorem}\label{thm:1} Let the basic assumptions on the system (\ref{eq:cl}) hold as stated in Section \ref{sec:1} and the exact solution be piecewise Lipschitz continuous with a finite number of isolated jump discontinuities. Let ${t\mapsto u_{\tau,h}(t,\cdot)\in\mathcal{D}}$ be an approximate solution to (\ref{eq:cl}), satisfying the weak consistency requirements  (\textbf{LC}), (\textbf{WC}), (\textbf{WES}), as well as the stability estimates in (\textbf{BV}). Then the following a priori error estimate holds,
	\begin{equation}
		\|u_{\tau,h}(T,\cdot)-\mathcal{S}_T(u_0)\|_{L^1(\R)}\lesssim h^{1/3}.
	\end{equation}
\end{theorem}

\begin{remark}
	In the presence of developing rarefaction waves and shock interactions in the exact solution, additional assumptions on the local oscillation decay of the numerical solution near developing rarefaction waves and the smallness of time intervals where shock interactions make us lose track of the shock curves must be made to generalize Theorem \ref{thm:1}. Under the heuristically appropriate scaling assumptions discussed in \cite{BressanChiriShen21}, cf.~Remark 4 there, an error bound of order $h^{1/3}\ln(h)$ can be derived. 
\end{remark}

\section{Finite Volume methods}\label{sec:fvmethod}
 Since we consider finite volume methods, the finite element space $V_h$ consists of piecewise constant functions on the space grid, $p(x)=\sum_j\mathbf{1}_{I_j}(x)p_j$ with $v_j\in\mathbb{R}^n$.

Given $u^n\in V_h$, the numerical solution $u^{n+1}\in V_h$ at time $t^{n+1}$ is given by
\begin{equation}\label{eq:scheme}
	\frac{u^{n+1}_j-u^n_j}{\tau} = - \frac{F^n_{j+1/2}(u^n
			_j,u^n_{j+1})-F^n_{j-1/2}(u^n_{j-1},u^n_j)}{h},
\end{equation}
where $F^n_{j+1/2}(a,b):=\widehat{F}(u^n_j,u^n_{j+1};a,b)$ is a numerical flux function such that the first two arguments of $\widehat{F}$ are understood as determining the numerical viscosity. Furthermore, we assume that the numerical solution assumes only values in a convex domain $\mathcal{O}\subset\R^n$.

To facilitate the proof of entropy stability, we will make use of the entropy variables: When $\eta$ is strictly convex on $\mathcal{O}$, there exists a convex domain $\widetilde{\mathcal{O}}\subset\R^n$ such that
\begin{equation}\label{eq:entrvar}
		v:\mathcal{O}\to\widetilde{\mathcal{O}},~u\mapsto \nabla\eta(u)
\end{equation}
is a diffeomorphism whose inverse we denote by $u:\widetilde{\mathcal{O}}\to\mathcal{O}$. With this, we may reformulate the numerical flux in terms of the entropy variables, i.e.
\begin{align}\label{eq:entrvar2}
	F^n_{j+1/2}(u^n_j,u^n_{j+1})=F^n_{j+1/2}(u(v^n_j),u(v^n_{j+1}))
	~\widehat{=}~\widetilde{F}^n_{j+1/2}(v^n_j,v^n_{j+1}),
\end{align}
where $v^n_j~\widehat{=}~v(u^n_j)$.

For the purpose of our analysis, the numerical flux shall satisfy the following conditions for all $j\in\mathbb{Z}$:
\begin{enumerate}[(i)]
	\item It is (first-order) \textit{consistent} with the physical flux $f$, i.e.~$F^n_{j+1/2}(p,p)=f(p)$ for all $p\in\mathcal{O}$..
	\item It is continuously differentiable and \textit{uniformly strictly monotone} in the entropy variables, i.e.~
	\begin{enumerate}[(1.)]
		\item The left-state Jacobian
			$D^l \widetilde{F}^n_{j-1/2}(~\cdot~,v^n_j)\left(=D\left[\widehat{F}\left(u^n_{j-1},u^n_j;u(\cdot),u^n_j\right)\right]\right)$
		is symmetric positive definite on the connecting line from $v^n_{j-1}$ to $v^n_j$ with uniformly bounded eigenvalues.
		\item The right-state Jacobian $D^r \widetilde{F}^n_{j+1/2}(v^n_j,~\cdot~)\left(=D\left[\widehat{F}\left(u^n_j,u^n_{j+1};u^n_j,u(\cdot)\right)\right]\right)$ is symmetric negative definite on the connecting line from $v^n_j$ to $v^n_{j+1}$ with uniformly bounded eigenvalues.
	\end{enumerate}
	\item It is \textit{second-order consistent}, i.e. there holds
	\begin{equation}\label{eq:fluxcons2}
		D^lF^n_{j-1/2}(~\cdot~,u^n_j)+D^rF^n_{j-1/2}(u^n_{j-1},~\cdot~)=Df(\cdot)
	\end{equation}
 along the curve $\left(u(v^n_{j-1}+s(v^n_j-v^n_{j-1}))\right)_{s\in[0,1]}$.
\end{enumerate}

\begin{remark}
	The most restrictive condition on the numerical flux is the strict monotonicity condition (ii). For one, this property might be violated due to the form of the physical flux function $f$, e.g.~in the case of vacuum states for the Euler system since eigenvalues blow up. On the other hand, if the eigenvalues of the numerical flux Jacobian are not bounded away from zero a subsequent entropy analysis might be impeded, thus numerical viscosity must be chosen large enough or an entropy fix has to be considered. 
\end{remark}
\begin{example}[HLL flux without entropy fix]\label{rem:HLL}
	 As a first example, consider the upwind HLL flux
	\begin{equation*}
		F^n_{j+1/2}(p,q)=\begin{cases}
			f(p), & s^- \geq 0\\
			\frac{s^+ f(p) - s^- f(q) - s^+s^- (p-q)}{s^+-s^-}, & s^- \leq 0 \leq s^+\\
			f(q), & s^+ \leq 0
		\end{cases},
	\end{equation*}
	where $s^-$ is the minimal and $s^+$ the maximal eigenvalue of the flux Jacobian $Df$ along the path $(u(v^n_{j,j+1}(s)))_{s\in[0,1]}$, where $v^n_{j,j+1}(s)=v^n_j+s(v^n_{j+1}-v^n_j)$. In the upwinding case $0=s^-<s^+$, the left-state Jacobian w.r.t.~entropy variables is just
	\begin{equation*}
		Df(u(v^n_{j,j+1}(s)))(\nabla^2\eta)^{-1}(v^n_{j,j+1}(s)),
	\end{equation*}
	which is symmetric due to (\ref{eq:symmetry}), but only positive semi-definite since $s^-=0$. It can be seen from the analysis in Section \ref{subsec:entrstab} that this might not yield enough entropy dissipation to counter anti-dissipative effects of a forward Euler time discretization. However, if the points where the minimal eigenvalue of $Df$ is between zero and $h^\alpha$, $\alpha\in(0,1)$, do not lie dense along the path $(u(v^n_{j,j+1}(s)))_{s\in[0,1]}$, entropy stability can be guaranteed under a CFL condition of the form $\tau\leq c h^{1+\alpha}$ where $c>0$ is -- among other dependencies -- inversely proportional to the "denseness" of those points in the worst case among all $j\in\mathbb{Z}$; similar considerations apply to the cases $s^-\leq 0\leq s^+$ and $s^+\leq 0$.
\end{example}

\begin{example}[HLL flux with entropy fix]
Let $\delta>0$ be some treshold. We define the HLL flux with entropy fix as follows.	
	\begin{equation*}
		F^n_{j+1/2}(p,q)=\begin{cases}
			f(p), & s^- \geq \delta\\
			\frac{(s^++\delta) f(p) + \delta f(q) + (s^++\delta)\delta (p-q)}{s^+ +2\delta}, & s^-\in[0,\delta)\\
			\frac{(s^++\delta) f(p) - (s^--\delta) f(q) - (s^++\delta)(s^--\delta) (p-q)}{s^+-s^-+2\delta}, & s^- \leq 0 \leq s^+\\
			\frac{\delta f(p) - (s^--\delta) f(q) - \delta(s^--\delta) (p-q)}{-s^-+2\delta}, &s^+\in(-\delta,0]\\
			f(q), & s^+ \leq -\delta
		\end{cases},
	\end{equation*}
	with $s^-,s^+$ as in Remark \ref{rem:HLL}. It can be easily verified that this numerical flux is uniformly strictly monotone as demanded with exception of the pure upwind cases, where either the left- or right-state Jacobian will vanish. However, the upwinding property suffices to prove entropy stability since here the propagation speeds $s^-,s^+$ are bounded away from zero by $\delta$.
\end{example}

\begin{remark}[Second-order consistency]
	The second-order consistency relation (\ref{eq:fluxcons2}) is satisfied by a variety of different numerical fluxes, e.g.~ Lax-Friedrichs/Rusanov-type, Enquist-Osher, as well as the HLL flux. In fact, consider the following (formal) equation that is due to mere consistency of a generic numerical flux function $F$,
	\begin{equation*}
		\int_{a}^{b}Df(s)~ds=\int_{a}^{b}D^lF(s,b)+D^rF(a,s)~ds,
	\end{equation*}
	which holds for all states $a,b$. If now $D^lF(s;b)$ does not depend on $b$ and $D^rF(a;s)$ does not depend on $a$, the second order consistency relation holds, since $a,b$ can be chosen arbitrarily. 
\end{remark}

\begin{example}
A Lax-Friedrichs/Rusanov-type flux has the form,
\begin{equation*}
	F^n_{j+1/2}(p,q)=\frac{f(p)+f(q)}{2}+\frac{1}{2}|A^n_{j+1/2}|(p-q),
\end{equation*}
with a stabilization matrix $A^n_{j+1/2}=A(u^n_j,u^n_{j+1})$; note that also the Enquist-Osher flux (with entropy fix) can be cast in this form. Monotonicity then depends on a suitable choice of the stabilization matrix and 
\begin{equation*}
	D^l F_{j+1/2}(p,q)=\frac{1}{2}(Df(p)+|A^n_{j+1/2}|),\quad D^r F_{j+1/2}(q,p)=\frac{1}{2}(Df(p)-|A^n_{j+1/2}|),
\end{equation*}
thus (\ref{eq:fluxcons2}) holds.	
\end{example}	

\begin{example}
	If we consider the HLL flux from above in the case $s^-\leq 0 \leq s^+$, then
	\begin{equation*}
		D^lF^n_{j+1/2}(p,q)=\frac{s^+(Df(p)-s^-\operatorname{Id})}{s^+-s^-},\quad D^rF^n_{j+1/2}(q,p)=\frac{s^-(s^+\operatorname{Id}-Df(p))}{s^+-s^-},
	\end{equation*}
	and adding the two expressions thus yields (\ref{eq:fluxcons2}).
\end{example}

\section{Analysis of Finite Volume methods}\label{sec:fvanalysis}
In this section, we will show that the finite volume solution,
\begin{equation}\label{eq:u}
	u(t,x):=u^n_j,\quad t\in[t^n,t^{n+1}),\quad x\in I_j,
\end{equation}
satisfies (\textbf{LC}), (\textbf{WC}) and (\textbf{WES}).

For shorthand notation we will abbreviate $u(t^n,\cdot)=u^n$ as well as $\varphi(t^n,\cdot)=\varphi^n$.

\subsection{Lipschitz continuity in time}
By definition of the numerical scheme (\ref{eq:scheme}) and Lipschitz continuity of the numerical flux, we have
\begin{align}
	\|u^{n+1}-u^n\|_{L^1(\R)}\lesssim\tau\sum_j|u^n_{j+1}-u^n_j|+|u^n_j-u^n_{j-1}|\lesssim \tau \totvar(u(t^n,\cdot)),
\end{align}
and thus (\textbf{LC}) follows inductively.

\subsection{Weak consistency}\label{sec:residual}
Let $\varphi\in C^1_c(\R^2)$. On the one hand, by definition of $u$ in (\ref{eq:u}), we have 
\begin{equation}\label{eq:WC1}
	\int_{t^n}^{t^{n+1}}\int_\R u\varphi_t~dx~dt - \int_\R u^{n+1}\varphi^{n+1}-u^n\varphi^n~dx=-\int_\R (u^{n+1}-u^n)\varphi^{n+1}~dx,
\end{equation}
and the numerical scheme (\ref{eq:scheme}) yields
\begin{equation*}
	-\int_\R (u^{n+1}-u^n)\varphi^{n+1}~dx=\tau\sum_j (F^n_{j+1/2}(u^n
	_j,u^n_{j+1})-F^n_{j-1/2}(u^n_{j-1},u^n_j))\langle\varphi^{n+1}\rangle_j,
\end{equation*}
where 
\begin{equation}
	\langle\varphi^{n+1}\rangle_j:=\frac{1}{h}\int_{I_j}\varphi^{n+1}~dx.
\end{equation}
On the other hand,
\begin{equation}\label{eq:WC2}
	\int_{t^n}^{t^{n+1}}\int_\R f(u)\varphi_x~dx = \tau\sum_j f(u^n_j)(\overline{\varphi}_{j+1/2}-\overline{\varphi}_{j-1/2}),
\end{equation}
where
\begin{equation}\label{eq:timeaverage}
	\overline{\varphi}_{j+1/2}:=\frac{1}{\tau}\int_{t^n}^{t^{n+1}}\varphi_{j+1/2}~dt.
\end{equation}
Note that
\begin{equation}\label{eq:approx}
	|\langle\varphi^{n+1}\rangle_j- \overline{\varphi}_{j-1/2}|\leq\frac{1}{2}(h\|\varphi_x\|_\infty+\tau\|\varphi_t\|_\infty).
\end{equation}
Now, similarly as in \cite{GiesselmannSikstel25}, by subtracting
\begin{equation*}
	0=\tau\sum_j F^n_{j+1/2}(u^n_j,u^n_{j+1})(\overline{\varphi}_{j+1/2}-\overline{\varphi}_{j-1/2}) + (F^n_{j+1/2}(u^n_j,u^n_{j+1})-F^n_{j-1/2}(u^n_{j-1},u^n_j))\overline{\varphi}_{j-1/2},
\end{equation*}
from (\ref{eq:WC1})+(\ref{eq:WC2}) we obtain
\begin{align}
	&\left|\int_{t^n}^{t^{n+1}}\int_\R u\varphi_t + f(u)\varphi_x~dx - \int_\R u^{n+1}\varphi^{n+1}-u^n\varphi^n~dx\right|\nonumber\\
	=~& \tau\sum_j \mathlarger{\mathlarger{|}}(F^n_{j+1/2}(u^n
	_j,u^n_{j+1})-F^n_{j-1/2}(u^n_{j-1},u^n_j))(\langle\varphi^{n+1}\rangle_j-\overline{\varphi}_{j-1/2})\nonumber\\
	&\qquad+ (f(u^n_j)-F^n_{j+1/2}(u^n
	_j,u^n_{j+1}))(\overline{\varphi}_{j+1/2}-\overline{\varphi}_{j-1/2})\mathlarger{\mathlarger{|}}\nonumber\\
	\lesssim~&\tau\sum_j \left(|u^n_{j+1}-u^n_j|+|u^n_j-u^n_{j-1}|\right)\left(h\|\varphi_x\|_\infty+\tau\|\varphi_t\|_\infty\right)+ h|u^n_{j+1}-u^n_j|\|\varphi_x\|_\infty\nonumber\\
	\lesssim~& h\|\varphi\|_{W^{1,\infty}(\R^2)}\cdot\tau\totvar(u(t^n,\cdot)),
\end{align}
due to consistency and Lipschitz continuity of the numerical flux as well as (\ref{eq:approx}) and a generic CFL condition of the form $\tau\lesssim h$.

\subsection{Weak entropy stability}\label{subsec:entrstab}
Let $\varphi\in C^1_c(\R^2)$ such that $\varphi\geq0$. We have on the one hand
\begin{align}
	&\int_{t^n}^{t^{n+1}}\int_\R \eta(u)\varphi_t~dx - \int_\R \eta(u^{n+1})\varphi^{n+1}-\eta(u^n)\varphi^n~dx\nonumber\\=&-\int_\R (\eta(u^{n+1})-\eta(u^n))\varphi^{n+1}~dx,
\end{align}
and the numerical scheme (\ref{eq:scheme}) together with the convexity of $\eta$ yield,
\begin{align}\label{eq:entropy scheme}
	-\int_\R (\eta(u^{n+1})-\eta(u^n))\varphi^{n+1}~dx\geq&~ \tau\sum_j (\nabla\eta(u^{n+1}_j)-\nabla\eta(u^n_j))\left(F(u^n_j,u^n_{j+1})-F(u^n_{j-1},u^n_j)\right)\langle\varphi^{n+1}\rangle_j\nonumber\\
	&\qquad+ \nabla\eta(u^n_j)\left(F(u^n_j,u^n_{j+1})-F(u^n_{j-1},u^n_j)\right)\langle\varphi^{n+1}\rangle_j,
\end{align} 
where we suppressed the subscripts of the numerical flux terms, e.g.~$F(u^n_j,u^n_{j+1})\hat{=}F^n_{j+1/2}(u^n_j,u^n_{j+1})$, for ease of notation.

\begin{remark}
	Note that 
	\begin{align}\label{eq:antidiss}
		(\nabla\eta(u^{n+1}_j)-\nabla\eta(u^n_j))\left(F(u^n_j,u^n_{j+1})-F(u^n_{j-1},u^n_j)\right)\leq 0,
	\end{align}
	due to the convexity of $\eta$. This term thus constitutes an anti-dissipative effect that originates from the forward Euler in time discretization. Since a bound from below as demanded by (\textbf{WES}) is not accessible for this term in isolation, it must be compensated by a dissipative mechanism of the numerical scheme.
\end{remark}

Let us now focus on the second term on the RHS of (\ref{eq:entropy scheme}). For one, we want to relate it to the exact entropy flux terms and, as mentioned in the remark above, balance it with the anti-dissipative effect of the explicit time discretization. For this purpose, we use a similar idea as in the multi-dimensional scalar setting considered in \cite{DednerMakridakisOhlberger07}, where we suitably add and subtract $f(u^n_j)$ and -- in contrast to the scalar case -- after moving to entropy variables, $v^n_j\hat{=}v(u^n_j)$ (cf.~(\ref{eq:entrvar})-(\ref{eq:entrvar2})), identify two contributions, the first one being
\begin{align}\label{eq:2dissip}
	\nabla\eta(u^n_j)\left(F(u^n_j,u^n_{j+1})-f(u^n_j)\right)=&~Q^r_{j+1/2}-q(u^n_j)\nonumber\\ &- \int_{0}^{1}s(v^n_{j+1}-v^n_j)^\intercal D^r \widetilde{F}(v^n_j,v^n_{j,j+1}(s))(v^n_{j+1}-v^n_j)~ds,
\end{align}
with $v^n_{j,j+1}(s):=v^n_j+s(v^n_{j+1}-v^n_j)$, and where
\begin{align}\label{eq:numentrflux2}
	Q^r_{j+1/2}:= &\int_0^1 v^n_{j,j+1}(s)^\intercal D^r\widetilde{F}(v^n_j,v^n_{j,j+1}(s))(v^n_{j+1}-v^n_j)~ds+ q(u^n_j),
\end{align}
is a consistent numerical entropy flux. Note that due to strict monotonicity of the numerical flux w.r.t.~entropy variables, the third term on the right hand side of (\ref{eq:2dissip}) is positive.

The second contribution is expressed similarly, as
\begin{align}\label{eq:1dissip}
	\nabla\eta(u^n_j)(f(u^n_j)-F(u^n_{j-1},u^n_j))=&~q(u^n_j)-Q^l_{j-1/2}\nonumber\\ &+ \int_{0}^{1}(1-s)(v^n_j-v^n_{j-1})^\intercal D^l \tilde{F}(v^n_{j-1,j}(s),v^n_j)(v^n_j-v^n_{j-1})~ds,
\end{align}
where
\begin{align}\label{eq:numentrflux1}
	Q^l_{j-1/2}:=&-\int_0^1 v^n_{j-1,j}(s)^\intercal D^l\widetilde{F}(v^n_{j-1,j}(s),v^n_j)(v^n_{j}-v^n_{j-1})~ds+ q(u^n_j),
\end{align}
is another consistent numerical entropy flux.
As before, due to monotonicity of the numerical flux w.r.t.~entropy variables, the third term on the RHS of (\ref{eq:1dissip}) is positive.

In order to inspect the difference between the two numerical entropy fluxes in (\ref{eq:numentrflux2}) and (\ref{eq:numentrflux1}), we note that
\begin{align}
	&v^n_{j-1,j}(s)^\intercal \left(D^l\widetilde{F}(v^n_{j-1,j}(s),v^n_j)+D^r\widetilde{F}(v^n_{j-1},v^n_{j-1,j}(s))\right)(v^n_{j}-v^n_{j-1})\nonumber\\
	=~&\nabla\eta(u(v^n_{j-1,j}(s)))Df(u(v^n_{j-1,j}(s)))(\nabla^2\eta)^{-1}(v^n_{j-1,j}(s))(v^n_{j}-v^n_{j-1})\nonumber\\
	=~&\frac{d}{ds} q(u(v^n_{j-1,j}(s))),
\end{align}
due to the second order consistency relation (\ref{eq:fluxcons2}), yielding
\begin{align}\label{eq:numentrfluxid}
	Q^r_{j-1/2}-Q^l_{j-1/2}= &\int_0^1 \frac{d}{ds} q(u(v^n_{j-1,j}(s)))~ds-\left(q(u^n_j)-q(u^n_{j-1})\right)=0,
\end{align}
We will thus drop the superscript for the numerical entropy flux.

It remains to balance the numerical entropy flux terms with the physical entropy flux terms, as well as the anti-dissipative term due to (\ref{eq:antidiss}) with the dissipative terms established in (\ref{eq:2dissip}) and (\refeq{eq:1dissip}).

The physical entropy flux terms are
\begin{align}
	\int_{t^n}^{t^{n+1}}\int_\R q(u)\varphi_x~dx =& \tau\sum_j q(u^n_j)(\overline{\varphi}_{j+1/2}-\overline{\varphi}_{j-1/2})\nonumber\\
	=& \tau\sum_j (q(u^n_j)-Q_{j+1/2})(\overline{\varphi}_{j+1/2}-\overline{\varphi}_{j-1/2})-(Q_{j+1/2}-Q_{j-1/2})\overline{\varphi}_{j-1/2},
\end{align}
where, similarly to Section \ref{sec:residual}, the RHS has the appropriate form to commence bounding the entropy residual, since
\begin{equation}\label{eq:QLipschitz}
	|q(u^n_j)-Q_{j+1/2}|\leq\lambda_{\max}(|Df|)\sup\limits_{u\in\mathcal{O}}|\nabla\eta(u)|\cdot|u^n_{j+1}-u^n_{j}|
\end{equation}

Lastly, we set off the anti-dissipative with the dissipative terms, e.g., on the one hand,
\begin{align}
	&(F(u^n_j,u^n_{j+1})-F(u^n_{j-1},u^n_j))\cdot(\nabla\eta(u^{n+1}_j)-\nabla\eta(u^n_j))\nonumber\\
	\geq& -2\frac{\tau}{h}\lambda_{\max}(\nabla^2\eta)\left(\lambda_{\max}(|D^r\widetilde{F}|)^2|v^n_{j+1}-v^n_j|^2+\lambda_{\max}(|D^l\widetilde{F}|)^2|v^n_j-v^n_{j-1}|^2\right)\nonumber\\
	\geq&-2\frac{\tau}{h}\frac{\lambda_{\max}(\nabla^2\eta)}{\lambda_{\min}(\nabla^2\eta)^2} \left(\lambda_{\max}(|D^rF|)^2|v^n_{j+1}-v^n_j|^2+\lambda_{\max}(|D^lF|)^2|v^n_j-v^n_{j-1}|^2\right)\label{eq:antidissipative},
\end{align}
and, on the other hand,
\begin{align}
	&\int_{0}^{1}-s(v^n_{j+1}-v^n_j)^\intercal D^r \tilde{F}(v^n_j,v^n_{j,j+1}(s))(v^n_{j+1}-v^n_j)~ds\nonumber\\
	\geq&~\frac{1}{2}\lambda_{\min}(|D^r\widetilde{F}|)|v^n_{j+1}-v^n_{j}|^2\geq \frac{\lambda_{\min}(|D^rF|)}{2\lambda_{\max}(\nabla^2\eta)}|v^n_{j+1}-v^n_{j}|^2,\label{eq:dissipative1}\\
	&\int_{0}^{1}(1-s)(v^n_{j}-v^n_{j-1})^\intercal D^l \tilde{F}(v^n_{j-1,j}(s),v^n_j)(v^n_{j}-v^n_{j-1})~ds\nonumber\\
	\geq&~\frac{1}{2}\lambda_{\min}(|D^l\widetilde{F}|)|v^n_{j}-v^n_{j-1}|^2\geq \frac{\lambda_{\min}(|D^lF|)}{2\lambda_{\max}(\nabla^2\eta)}|v^n_{j+1}-v^n_{j}|^2\label{eq:dissipative2},
\end{align}
from which we can read off a net entropy dissipative result as long as the following CFL condition holds,
\begin{equation}\label{eq:cfl}
	\frac{\tau}{h}\leq \frac{1}{4}\frac{\lambda_{\min}(\nabla^2\eta)^2}{\lambda_{\max}(\nabla^2\eta)^2}\cdot\frac{\min(\lambda_{\min}(|D^lF|),\lambda_{\min}(|D^rF|))}{\max(\lambda_{\max}(|D^lF|),\lambda_{\max}(|D^rF|))^2}.
\end{equation}
We thus conclude this section with
\begin{lemma}
	A finite volume solution, (\ref{eq:u}), computed from (\ref{eq:scheme}) with a first- and second-order consistent numerical flux that is strictly monotone satisfies (\textbf{LC}), (\textbf{WC}), (\textbf{WES}) under the CFL condition (\ref{eq:cfl}).
\end{lemma}

\begin{cor}\label{cor:fv}
	If a finite volume solution, (\ref{eq:u}), computed from (\ref{eq:scheme}) with a first- and second-order consistent numerical flux that is strictly monotone, fulfills the stability estimates in (\textbf{BV}), then the a priori error estimate,
	\begin{equation}
		\|u(T,\cdot)-S_T(u_0)\|_{L^1(\R)}\lesssim h^{1/3},
	\end{equation}
	to the exact solution to (\ref{eq:cl}), as specified in Theorem \ref{thm:1}, holds true under the CFL condition (\ref{eq:cfl}).
\end{cor}
\begin{remark}
	 In \cite{Tadmor03} a comparable time step restriction w.r.t.~the (cube of the) condition number of the entropy Hessian is formulated to ensure entropy stability.
\end{remark}

\begin{remark}
	Regarding the behaviour of the constant $\alpha:=\frac{\lambda_{\min}(\nabla^2\eta)^2}{\lambda_{\max}(\nabla^2\eta)^2}$, take for example the isothermal Euler system. For this system, $\alpha \ll 1$ near the vacuum case or for velocities $\gg1$. 
\end{remark}
\begin{remark}
	Regarding the constant $\beta:=\frac{\min(\lambda_{\min}(|D_1F|),\lambda_{\min}(|D_2F|))}{\max(\lambda_{\max}(|D_1F|),\lambda_{\max}(|D_2F|))^2}$, take for example a numerical flux of the form
	\begin{equation*}
		F_{j+1/2}(u,v)=\frac{f(u)+f(v)}{2}+\frac{\mu}{2}|A_{j+1/2}|(u-v),
	\end{equation*} 
with $\mu>1$ and $|A_{j+1/2}|=\lambda_{\max}(|Df|)\operatorname{Id}^{n\times n}$, where the maximal eigenvalue of $|Df|$ is chosen along the path $(u(v^n_{j,j+1}(s)))_{s\in[0,1]}$. Then
\begin{equation*}
	\beta=2\frac{\mu-1}{(\mu+1)^2}\cdot\frac{1}{\lambda_{\max}(|Df|)},
\end{equation*}
which is maximized for $\mu=3$.
\end{remark}

\section{Forward Euler $P^k$-DG methods without limiting}\label{sec:without limiting}
In this section, we show that a single forward Euler time step for a $P^k$-DG method respects the inequalities in (\textbf{LC}), (\textbf{WC}) and (\textbf{WES}). However, since forward Euler DG methods are known to be linearly unstable, the considerations should be regarded as preparation for Section \ref{sec:with limiting}.

In the following, $\langle\cdot,\cdot\rangle$ denotes the scalar product in $L^2(\R;\R^n)$ and $\llbracket p\rrbracket_{j+1/2}$ denotes the jump $p(x_{j+1/2}^+)-p(x_{j+1/2}^-)$ of a function at the grid point $x_{j+1/2}$. Let $V_h^k$ be the space of piecewise polynomials of order $k$ on the space grid. Given $u^n\in V_h^k$, find $u^{n+1}\in V_h^k$, such that

\begin{align}
	\left\langle \frac{u^{n+1} - u^n}{\tau}, p\right\rangle =& ~\langle f(u^n), p_x\rangle	+ \sum_j \widehat{f(u^n)}_{j+1/2}\llbracket p\rrbracket_{j+1/2}\nonumber\\
	=& - \langle f(u^n)_x, p\rangle-\sum_j(\widehat{f(u^n)}_{j+1/2}-f(u^n)_{j+1/2}^-)p_{j+1/2}^-\nonumber\\
	& +\sum_j(\widehat{f(u^n)}_{j-1/2}-f(u^n)_{j-1/2}^+)p_{j-1/2}^+,
\end{align}
for all $p\in V_h^k$. Here, $\widehat{f(u^n)}_{j+1/2}:=F\left(u^n(x_{j+1/2}^-),u^n(x_{j+1/2}^+);u^n(x_{j+1/2}^-),u^n(x_{j+1/2}^+)\right)$ denotes the numerical flux, where $F$ is a numerical flux function as given in Section \ref{sec:fvmethod}.

For the analysis of the scheme we set,
\begin{equation}
	u(t,\cdot):=u^n,\quad\forall t\in[t^n,t^{n+1}).
\end{equation}

We will only discuss entropy stability of one forward Euler step in detail since consistency follows from similar considerations.

\subsection{Entropy production}\label{subsec:entropy production}

Let $\varphi\in C^1_c(\R^2)$ such that $\varphi\geq0$. When trying to compute the entropy violation for the DG scheme, of course one has to introduce projection errors since the entropy gradient will not live in the DG space in general. Thus, we have to investigate the following contributions, where we may drop the superscript for $\varphi(=\varphi^{n+1})$ for notational convenience:
\begin{align}\label{eq:entropysplitting}
	-\left\langle\nabla\eta(u^n)\varphi,u^{n+1}-u^n\right\rangle =& -\left\langle\nabla\eta(u^n)(\varphi-\Pi^0_h\varphi),u^{n+1}-u^n\right\rangle\nonumber\\
	&-\left\langle\Pi^k_h[\nabla\eta(u^n)]\Pi^0_h\varphi,u^{n+1}-u^n\right\rangle,
\end{align}
where $\Pi^k_h$ is the $L^2$-projection into $V^k_h$.
\begin{remark}\label{remark:phiregularity}
	In \cite{BressanChiriShen21}, the testfunctions $\varphi$ that are constructed to yield an estimate for the numerical error are merely $W^{1,\infty}(\R^2)$ and also no grid function, in general. As such, it makes no sense to use higher order interpolation for $\varphi$ in (\ref{eq:entropysplitting}).
\end{remark}

Let us first study the projection error term, i.e.
\begin{align}
	\mathcal{E}^0 :=|\left\langle\nabla\eta(u^n)(\varphi-\Pi^0_h\varphi),u^{n+1}-u^n\right\rangle|\leq&~ \sup\limits_{u\in\mathcal{O}}|\nabla\eta(u)|\sum_j \|\varphi-\Pi^0_h\varphi\|_{L^2(I_j)}\|u^{n+1}-u^n\|_{L^2(I_j)}\nonumber\\
	\leq &~ \sup\limits_{u\in\mathcal{O}}|\nabla\eta(u)|\cdot \|\varphi_x\|_\infty\sum_j h^{3/2}\|u^{n+1}-u^n\|_{L^2(I_j)}.
\end{align}
Now observe that in fact, due to an inverse inequality and $L^\infty-L^2$-estimate for discrete functions, we have
\begin{align}\label{eq:timestep difference L^2}
	\|u^{n+1}-u^n\|_{L^2(I_j)} \lesssim&~ \tau h^{-1/2} (\|Df\|_\infty +\|D\widehat{f}\|_\infty)(\|u^n_x\|_{L^1(I_j)}+|\llbracket u^n\rrbracket|_{j-1/2}+|\llbracket u^n\rrbracket|_{j+1/2}),
\end{align}
thus
\begin{equation}\label{eq:E^0}
	\mathcal{E}^0\lesssim \tau h \|\varphi_x\|_\infty \totvar(u^n).
\end{equation}

\begin{remark}
	The estimate (\ref{eq:timestep difference L^2}) also yields the (\textbf{LC}) property.
\end{remark}

\begin{remark}
	Let us note that the error bound (\ref{eq:E^0}) does not contain the $L^\infty$-norm of $u^n_x$ as a prefactor.
\end{remark}

For the investigation of the discrete entropy evolution, setting $\varphi_j:=\Pi^0_h\varphi|_{I_j}\in\R$, we write
\begin{align}
	-\left\langle\Pi^k_h[\nabla\eta(u^n)]\Pi^0_h\varphi,u^{n+1}-u^n\right\rangle = -\sum_j \left\langle\Pi^k_h[\nabla\eta(u^n)]|_{I_j},u^{n+1}-u^n\right\rangle\varphi_j,
\end{align}
and thus obtain three contributions to study for each cell $I_j$, i.e.
\begin{align}
	-\left\langle\Pi^k_h[\nabla\eta(u^n)]|_{I_j},u^{n+1}-u^n\right\rangle= &~\tau \langle f(u^n)_x, \Pi^k_h[\nabla\eta(u^n)]|_{I_j}\rangle\nonumber\\
	&+\tau\Pi^k_h[\nabla\eta(u^n)]_{j+1/2}^-\left(\widehat{f(u^n)}_{j+1/2}-f(u^n)_{j+1/2}^-\right)\nonumber\\
	&+\tau\Pi^k_h[\nabla\eta(u^n)]_{j-1/2}^+\left(f(u^n)_{j-1/2}^+ -\widehat{f(u^n)}_{j-1/2}\right)\nonumber\\
	\hat{=}&~ T_j^1+T_j^2+T_j^3.
\end{align}

\subsubsection{Bounding $T^1_j$}
The first term we treat as follows,
\begin{align}
	T_j^1=& ~\tau \langle f(u^n)_x, \Pi^k_h[\nabla\eta(u^n)]|_{I_j}-\nabla\eta(u^n)|_{I_j}\rangle+\tau \langle f(u^n)_x, \nabla\eta(u^n)|_{I_j}\rangle,
\end{align}
while the projection error term can be bounded as,
\begin{align}
	\mathcal{E}_j^1:=&~\tau |\langle f(u^n)_x, \Pi^k_h[\nabla\eta(u^n)]|_{I_j}-\nabla\eta(u^n)|_{I_j}\rangle|\nonumber\\ 
	\lesssim&~\tau\|f'\|_\infty h^{1/2}\|u^n_x\|_{L^1(I_j)}|\nabla\eta(u^n)|_{H^1(I_j)},
\end{align}
thus we obtain
\begin{equation}\label{eq:E_j^1}
	\mathcal{E}_j^1\lesssim \tau h \|u^n_x\|_{L^\infty(I_j)}\|u^n_x\|_{L^1(I_j)}.
\end{equation}
On the other hand, due to the entropy compatibility condition (\ref{eq:entropy compatibility}), we have
\begin{equation}
	\langle f(u^n)_x, \nabla\eta(u^n)|_{I_j}\rangle = q(u^n)_{j+1/2}^- - q(u^n)_{j-1/2}^+,
\end{equation}
and thus in total,
\begin{equation}\label{eq:DG entropy production volume}
	 \sum_j T_j^1\varphi_j \gtrsim -\tau h\|\varphi\|_\infty\|u^n_x\|_\infty \totvar(u^n) + \tau \sum_j \left(q(u^n)_{j+1/2}^- - q(u^n)_{j-1/2}^+\right)\varphi_j.
\end{equation}

\subsubsection{Bounding $T^2_j$ and $T^3_j$} Let us now focus on the term $T^2_j$, for $T^3_j$ is handled analogously. We write
\begin{align}
	T^2_j =& ~\tau\left(\Pi^k_h[\nabla\eta(u^n)]_{j+1/2}^- -\nabla\eta(u^n)_{j+1/2}^-\right)\left(\widehat{f(u^n)}_{j+1/2}-f(u^n)_{j+1/2}^-\right)\nonumber\\
	&+\tau \nabla\eta(u^n)_{j+1/2}^-\left(\widehat{f(u^n)}_{j+1/2}-f(u^n)_{j+1/2}^-\right),
\end{align}
and bound the projection error term as
\begin{align}\label{eq:E^2_j}
	\mathcal{E}_j^2 :=& ~\tau\left|\left(\Pi^k_h[\nabla\eta(u^n)]_{j+1/2}^- -\nabla\eta(u^n)_{j+1/2}^-\right)\right|\left|\left(\widehat{f(u^n)}_{j+1/2}-f(u^n)_{j+1/2}^-\right)\right|\nonumber\\
	\lesssim&~ \tau h^{-1/2}\|\Pi^k_h[\nabla\eta(u^n)]-\nabla\eta(u^n)\|_{L^2(I_j)}\|D\widehat{f}\|_\infty |\llbracket u^n\rrbracket|_{j+1/2}\nonumber\\
	& + \tau h^{-1/2}\|\nabla\eta(u^n)-\widetilde{\Pi}_h[\nabla\eta(u^n)]\|_{L^2(I_j)}\|D\widehat{f}\|_\infty |\llbracket u^n\rrbracket|_{j+1/2}\nonumber\\
	\lesssim& ~\tau h \|u^n_x\|_{L^\infty(I_j)}|\llbracket u^n\rrbracket|_{j+1/2},
\end{align}
where we have employed a discrete trace estimate with the linear Gauss-Radau interpolant $\widetilde{\Pi}_hv|_{I_j}$ of $v$ that preserves the function value on the right boundary point, i.e. $\widetilde{\Pi}_hv|_{I_j}(x_{j+1/2})=v(x_{j+1/2})^-$; cf.~\cite{LeottaGiesselmann24}.

Summarizing, we thus have
\begin{align}\label{eq:DG entropy production boundary}
	\sum_j (T^2_j+T^3_j)\varphi_j \gtrsim& -\tau h\|\varphi\|_\infty\|u^n_x\|_\infty\totvar(u^n)\nonumber\\
	&+\tau\sum_j \nabla\eta(u^n)_{j+1/2}^-\left(\widehat{f(u^n)}_{j+1/2}-f(u^n)_{j+1/2}^-\right)\varphi_j\nonumber\\
	&+\tau \sum_j\nabla\eta(u^n)_{j-1/2}^+\left(f(u^n)_{j-1/2}^+ -\widehat{f(u^n)}_{j-1/2}\right)\varphi_j,
\end{align}
meaning that relevant entropy production is solely due to the inter-element flux terms as in the finite volume setting.

\subsection{Entropy production vs.~Entropy flux}

Up to suitably bounded projection errors, we have identified the following entropy production terms in Section \ref{subsec:entropy production}, cf. (\ref{eq:DG entropy production volume}) and (\ref{eq:DG entropy production boundary}), associated to each cell $I_j$,
\begin{align}
	&\tau \left[q(u^n)_{j+1/2}^- + \nabla\eta(u^n)_{j+1/2}^-\left(\widehat{f(u^n)}_{j+1/2}-f(u^n)_{j+1/2}^-\right)\right]\varphi_j^{n+1}\label{eq:j+1/2}\\
	&\tau \left[\nabla\eta(u^n)_{j-1/2}^+\left(f(u^n)_{j-1/2}^+ -\widehat{f(u^n)}_{j-1/2}\right)-q(u^n)_{j-1/2}^+\right]\varphi_j^{n+1}\label{eq:j-1/2}.
\end{align}
From the finite volume analysis in subsection \ref{subsec:entrstab}, we know that the sum of (\ref{eq:j+1/2}) and (\ref{eq:j-1/2}) corresponds to a numerical entropy flux (that is consistent and yields a Lipschitz estimate, cf.~(\ref{eq:QLipschitz})) difference $Q_{j+1/2}-Q_{j-1/2}$ plus some dissipative terms $D_{j+1/2},~D_{j-1/2}$; cf.~Section \ref{subsec:dissipation} for a brief discussion of the balance between anti-dissipative and dissipative terms. Since now
\begin{align}
	\int_{t^n}^{t^{n+1}}\int_\R q(u)\varphi_x~dx~dt =& \tau\sum_j -\langle q(u^n)_x,\overline{\varphi}-\Pi^0_h\overline{\varphi}\rangle_{I_j}\nonumber\\
	& +\tau\sum_j q(u^n)_{j+1/2}^-(\overline{\varphi}_{j+1/2}-\overline{\varphi}_j) - q(u^n)_{j-1/2}^+(\overline{\varphi}_{j-1/2}-\overline{\varphi}_j),
\end{align}
where $\overline{\varphi}_{j\pm1/2}$ is the time average of $\varphi(\cdot,x_{j\pm 1/2})$, cf.~(\ref{eq:timeaverage}), and 
\begin{align}
	-\langle q(u^n)_x,\overline{\varphi}-\Pi^0_h\overline{\varphi}\rangle_{I_j}\geq - h\|\varphi_x\|_\infty\sup\limits_{u\in\mathcal{O}}|\nabla q(u)| \cdot\|u^n_x\|_{L^1(I_j)},
\end{align}
while we can write
\begin{align}
	 q(u^n)_{j+1/2}^-(\overline{\varphi}_{j+1/2}-\overline{\varphi}_j)=&~(q(u^n)_{j+1/2}^- -Q_{j+1/2})(\overline{\varphi}_{j+1/2}-\overline{\varphi}_j)\nonumber\\ &+ Q_{j+1/2}(\overline{\varphi}_{j+1/2}-\overline{\varphi}_j),\\ 
	 - q(u^n)_{j-1/2}^+(\overline{\varphi}_{j-1/2}-\overline{\varphi}_j)=&-(q(u^n)_{j-1/2}^+ -Q_{j-1/2})(\overline{\varphi}_{j-1/2}-\overline{\varphi}_j)\nonumber\\ &- Q_{j-1/2}(\overline{\varphi}_{j-1/2}-\overline{\varphi}_j),
\end{align}
we see that
\begin{align}
	\int_{t^n}^{t^{n+1}}\int_\R q(u)\varphi_x~dx~dt \gtrsim& -\tau h \|\varphi_x\|_\infty\totvar(u^n)\nonumber\\
	 &- \tau\sum_j (Q_{j+1/2}-Q_{j-1/2})\overline{\varphi}_j,
\end{align}
which is clearly compatible with the numerical entropy flux difference that is obtained from summing (\ref{eq:j+1/2}) and (\ref{eq:j-1/2}).

\subsection{Anti-dissipative terms}\label{subsec:dissipation}
While the dissipative terms resulting from monotonicity properties of the numerical flux are essentially the same as in the finite volume setting, we want to briefly address the shape of the anti-dissipative terms due to forward Euler time integration, i.e.
\begin{align}
	-\left\langle\nabla\eta(u^{n+1})-\nabla\eta(u^n),(u^{n+1}-u^n)\varphi\right\rangle_{I_j}.
\end{align}
Of course, as done in (\ref{eq:entropysplitting}), we can work with $\Pi^0_h\varphi$ instead of $\varphi$. Then, by testing the numerical scheme with the $L^2$-projection of $\nabla\eta(u^{n+1})-\nabla\eta(u^n)$ into the DG space, we see that the resulting volume terms can be bounded using (\ref{eq:timestep difference L^2}), i.e.
\begin{align}
	\tau\langle f(u^n)_x,\Pi^k_h\left[\nabla\eta(u^{n+1})-\nabla\eta(u^n)\right]\rangle_{I_j} \geq& -\tau h^{1/2}\|u^n_x\|_{L^\infty(I_j)} \|\Pi^k_h\left[\nabla\eta(u^{n+1})-\nabla\eta(u^n)\right]\|_{L^2(I_j)}\nonumber\\
	\gtrsim& -\tau^2 \|u^n_x\|_{L^\infty(I_j)}\left(\|u^n_x\|_{L^1(I_j)}+|\llbracket u^n\rrbracket|_{j-1/2}+|\llbracket u^n\rrbracket|_{j+1/2}\right),
\end{align}
such that using a generic CFL condition of the form $\tau\lesssim h$ shows that the relevant anti-dissipation stems solely from the inter-element boundary flux terms as in the finite volume case. This closes the discussion of the forward Euler DG method.

\section{A two-stage RK $P^k$-DG method with limiting}\label{sec:with limiting}
Since the forward Euler DG method is known to be linearly unstable, we consider a two-stage RK $P^k$-DG method with limiting in this section. We will be able to recycle the discussion from Section \ref{sec:without limiting}. 

With space and time grid as before, given a datum $u^n\in V^k_h$, we seek $u^{n+1}\in V^k_h$, such that

\begin{align}\label{eq:P^klimitedscheme}
	\langle u^{n,1},p\rangle =&\langle \widetilde{u}^n,p\rangle + \tau\langle f(\widetilde{u}^n),p_x\rangle + \tau\sum_j \widehat{f(\widetilde{u}^n)}_{j+1/2}\llbracket p\rrbracket_{j+1/2}\\
	\langle u^{n+1},p\rangle =&\frac{1}{2}\langle \widetilde{u}^n+\widetilde{u}^{n,1},p\rangle + \frac{1}{2}\tau\langle f(\widetilde{u}^{n,1}),p_x\rangle + \frac{1}{2}\tau\sum_j \widehat{f(\widetilde{u}^{n,1})}_{j+1/2}\llbracket p\rrbracket_{j+1/2},
\end{align}
for all $p\in V^k_h$, where $\widetilde{v}\in V^k_h$ denotes the slope-limited version of $v\in V^k_h$ that is assumed to conserve the mean, i.e., there holds
\begin{align}\label{eq:limiter1}
	\Pi^0_h(v)=\Pi^0_h(\widetilde{v}),
\end{align}
and a natural assumption on the limited slopes,
\begin{equation}\label{eq:limiter2}
	\|\widetilde{v}_x\|_{L^\infty(I_j)}\lesssim \|v_x\|_{L^\infty(I_j)}.
\end{equation}

As before, e.g.~$\widehat{f(\widetilde{u}^n)}_{j+1/2}:=F\left(\widetilde{u}^n(x_{j+1/2}^-),\widetilde{u}^n(x_{j+1/2}^+);\widetilde{u}^n(x_{j+1/2}^-),\widetilde{u}^n(x_{j+1/2}^+)\right)$ denotes the numerical flux, where $F$ is a numerical flux function as given in Section \ref{sec:fvmethod}.

We set
\begin{equation}\label{eq:P^1-RKDG}
	u(t,\cdot):=\widetilde{u}^n,\quad t\in[t^n,t^{n+1}).
\end{equation}
 
Again, we only discuss entropy stability of the method since consistency follows from similar considerations.

\subsection{Entropy stability}
Let $\varphi\in C^1_c(\R^2)$ such that $\varphi\geq0$. Due to convexity, we have
\begin{align}
	\int_{t^n}^{t^{n+1}}\int_\R\eta(u)\varphi_t~dx~dt - \int_\R [\eta(u)\varphi]^{n+1}_n~dx\geq& -\int_\R \left(\eta(u^{n+1})-\eta((\widetilde{u}^n+\widetilde{u}^{n,1})/2)\right)\varphi^{n+1}~dx + I^1\nonumber\\
	& - \frac{1}{2}\int_\R (\eta(u^{n,1})-\eta(\widetilde{u}^n))\varphi^{n+1}~dx + I^2
\end{align}
with the perturbations due to limiting,
\begin{align}
	I^1:=&-\int_\R \left(\eta(\widetilde{u}^{n+1})-\eta(u^{n+1})\right)\varphi^{n+1}~dx\\
	I^2:=&-\frac{1}{2}\int_\R (\eta(\widetilde{u}^{n,1})-\eta(u^{n,1}))\varphi^{n+1}~dx.
\end{align}
Working in a cell $I_j$, denoting $u_j:=\Pi^0_h(u)|_{I_j}$ and dropping superscripts, Taylor expansion yields
\begin{align}\label{eq:sigma}
	\eta(u)-\eta(\widetilde{u})\geq &~\nabla\eta(u_j)(u-u_j)-\nabla\eta(u_j)(\widetilde{u}-u_j)\nonumber\\
	&+\frac{1}{2}\left(\lambda_{\min}(\nabla^2\eta)|u-u_j|^2-\lambda_{\max}(\nabla^2\eta)|\widetilde{u}-u_j|^2\right).
\end{align}
Now, due to orthogonality and the limiter properties (\ref{eq:limiter1})-(\ref{eq:limiter2}), we have
\begin{equation}
	\int_{I_j}(\nabla\eta(u_j)(u-u_j)-\nabla\eta(u_j)(\widetilde{u}-u_j))\varphi~dx\gtrsim -h^2\|\varphi_x\|_\infty \|u_x\|_{L^1(I_j)}.
\end{equation}
On the other hand, again using the limiter properties,
\begin{equation}\label{eq:suboptimal}
	\int_{I_j}\left(\lambda_{\min}(\nabla^2\eta)|u-u_j|^2-\lambda_{\max}(\nabla^2\eta)|\widetilde{u}-u_j|^2\right)\varphi~dx\gtrsim -h^2\|\varphi\|_\infty \|u_x\|_\infty \|u_x\|_{L^1(I_j)}.
\end{equation}

\begin{remark}
	In the scalar case, a change of variables shows that the difference of the quadratic terms resulting from Taylor expansion in (\ref{eq:sigma}) is in fact non-negative. It is not clear how to generalize this to the case of nonlinear systems.
\end{remark}

Summarizing, we thus arrive at
\begin{align}
	I^1\gtrsim&-h^2\|\varphi\|_{W^{1,\infty}}(1+\|u^{n+1}_x\|_\infty)\totvar(u^{n+1})\nonumber\\
	\gtrsim&-\tau h\|\varphi\|_{W^{1,\infty}}(1+\|\widetilde{u}^n_x\|_\infty) \totvar(\widetilde{u}^n),
\end{align}
where we have used a generic CFL condition of the form $\tau\lesssim h$ to bound local variations of $u^{n+1}$ in terms of local variations of $\widetilde{u}^n$. 

Similarly, one arrives at
\begin{align}
	I^2\gtrsim -\tau h\|\varphi\|_{W^{1,\infty}}(1+\|\widetilde{u}^n_x\|_\infty) \totvar(\widetilde{u}^n).
\end{align}

The essentials of the remaining work are to carry out the analysis of the resulting forward Euler steps for the two Runge-Kutta stages which we have done in Section \ref{sec:without limiting}.

\begin{lemma}
	A two-stage RK $P^k$-DG solution with abstract limiting, (\ref{eq:P^1-RKDG}), computed from (\ref{eq:P^klimitedscheme}) with a first- and second-order consistent numerical flux that is strictly monotone satisfies (\textbf{LC}), (\textbf{WC}), (\textbf{WES}) under the CFL condition (\ref{eq:cfl}).
\end{lemma}

\begin{cor}\label{cor:DG}
	If a two-stage RK $P^k$-DG solution with abstract limiting, (\ref{eq:P^1-RKDG}), computed from (\ref{eq:P^klimitedscheme}) with a first- and second-order consistent numerical flux that is strictly monotone, fulfills the stability estimates in (\textbf{BV}) and is such that the slopes remain bounded in $L^\infty$, uniformly in $\tau$ and $h$, then the a priori error estimate,
	\begin{equation}
		\|u(T,\cdot)-S_T(u_0)\|_{L^1(\R)}\lesssim h^{1/3},
	\end{equation}
	to the exact solution to (\ref{eq:cl}), as specified in Theorem \ref{thm:1}, holds true under the CFL condition (\ref{eq:cfl}).
\end{cor}

\bibliography{references}{}
\bibliographystyle{plain}

\end{document}